\documentclass[11pt]{article}
\usepackage{amscd,amsmath, amssymb, fancyhdr, color}

\newcommand{\version}{version 1.0}

\newcommand{\f}{\varphi}

\numberwithin{equation}{section}

\def\eqref#1{(\ref{#1})}

\newcommand{\C}{{\mathbb C}}
\newcommand{\R}{{\mathbb R}}
\newcommand{\Q}{{\mathbb Q}}
\renewcommand{\H}{{\mathbb H}}

\def\1{\sqrt{-1}\:}

\newcommand{\cntrct}                
{\hspace{2pt}\raisebox{1pt}{\text{$\lrcorner$}}\hspace{2pt}}

\renewcommand{\phi}{\varphi}
\renewcommand{\epsilon}{\varepsilon}
\renewcommand{\geq}{\geqslant}
\renewcommand{\leq}{\leqslant}

\newcounter{Mycounter}[section]
\newcounter{lemma}[section]
\renewcommand{\thelemma}{{Lemma \thesection.\arabic{lemma}}}
\newcommand{\lemma}{%
     \setcounter{lemma}{\value{Mycounter}}
     \refstepcounter{lemma}
     \stepcounter{Mycounter}
     {\noindent \bf \thelemma:\ }}

\newcounter{claim}[section]
\newcounter{sublemma}[section]
\newcounter{corollary}[section]
\newcounter{theorem}[section]
\renewcommand{\thetheorem}{{Theorem \thesection.\arabic{theorem}}}
\newcommand{\theorem}{%
     \setcounter{theorem}{\value{Mycounter}}
     \refstepcounter{theorem}
     \stepcounter{Mycounter}
     {\noindent \bf \thetheorem:\ }}

\newcounter{conjecture}[section]
\newcounter{proposition}[section]
\newcounter{definition}[section]
\newcounter{example}[section]
\newcounter{remark}[section]
\newcounter{problem}[section]
\newcounter{question}[section]
\makeatletter

\@addtoreset{equation}{section}
\@addtoreset{footnote}{section}
\makeatother

\def\O{\mathcal O}

\begin{document}
\begin{center}
{\LARGE\bf
 LCK metrics on Oeljeklaus-Toma manifolds vs Kronecker's theorem}\\[3mm]

Victor Vuletescu\footnote{Partially supported by CNCS –UEFISCDI, project
number PN-II-ID-PCE-2011-3-0118
\\[2mm]

{\bf Keywords:} Locally
conformally K\"ahler manifold,
number field, units.

{\bf 2000 Mathematics Subject
Classification:} { 53C55.}}
\\[4mm]

\end{center}

{\small
\hspace{0.15\linewidth}
\begin{minipage}[t]{0.7\linewidth}
{\bf Abstract} \\ A locally conformally K\"ahler (LCK) manifold is a manifold which is covered by a K\"ahler manifold, with the deck transform group acting by homotheties. We show that the search for LCK metrics on Oeljeklaus-Toma manifolds leads to a (yet another) variation  on  Kronecker's theorem on units. In turn, this implies that on Oeljeklaus-Toma manifold  associated to number fields with $2t$ complex embeddings and $s$ real embeddings with $s<t$ there is no LCK metric.
\end{minipage}
}

\section{Introduction}

\subsection{Locally conformally K\"ahler structures}

A {\bf locally conformally K\"ahler} (LCK) manifold is a 
complex manifold $X$, $\dim_\C X >1$, admitting
a K\"ahler covering $(\tilde X,  \tilde \omega)$, with
the deck transform group acting on $(\tilde X, \tilde \omega)$
by holomorphic homotheties. In other words, for all $\gamma \in \pi_1(X)\subset Aut(\tilde X)$ there exists some $\chi(\gamma)\in \R_{>0}$ such that 
$$\gamma^*(\tilde \omega)=\chi(\gamma)\tilde \omega.$$
The positive numbers $\chi(\gamma)$ are called {\em the automorphy factors} of $X.$

LCK manifolds were introduced in the late 70's by I Vaisman, in an attempt to exhibit interesting metrics on non-K\"ahler manifolds. Basically, Vaisman noticed that the fundamental group of a standard Hopf manifold $X$ (for simplicity, generated by $(z\mapsto 2z)$) acts on the standard flat metric $\omega_0$ on $\C^n\setminus \{0\}$ by homotheties; consequently, the metric $\frac{1}{\vert z\vert}\omega_0$ descends to the quotient. Thus, even if $X$ has no K\"ahler metric (for instance, since it has first Betti number equal to $1$), it still carries a  interesting metric, as  $\frac{1}{\vert z\vert}\omega_0$ is loccally conformal  to a K\"ahler one.

Deciding whether a given (compact) complex manifold belongs or not to the class of manifolds carrying  an LCK metric is a rather tricky problem. No general procedures can apply; this class is known not be closed under (even small!) deformations, and is still an open problem whether is closed or not under taking products or  finite quotients. On the other hand, no general (e.g. topological) restrictions are known; except for the non-simply-connectedness, only some mild restrictions are known on the fundamental group of a compact compact manifold that prevent it from having LCK metrics with additional properties (see e.g. \cite{OrVeFG}). 

Despite these difficulties, along the years, a rather suprising result e\-merged: almost all compact complex non-K\'ahler surfaces have LCK metrics! A rough chronological list would include (apart from standard Hopf surfaces from Vaisman's original paper): one class of Inoue surfaces (Tricerri, 1982, \cite{Tri}), general Hopf surfaces (Gauduchon-Ornea, 1998, \cite{GaOr}), elliptic surfaces and another class of Inoue surfaces (Belgun, 2000, \cite{Bel}) and eventually the only known examples of surfaces in Kodaira's class $VII_b$ with $b>0$, namely Kato surfaces (Brunella, 2010-2011, \cite{Bru1}, \cite{Bru2}). Let us mention, that the only class of non-K\"ahler surfaces  known so far not to admit LCK structures is a third class of Inoue surfaces (Belgun, \cite{Bel}) and, possibly, some hypothetical non-Kato surfaces in class $VII_b, b>0$ - which are also supposed,  by the {\em global spherical shells conjecture}, not to exist!

In higher dimensions, the only know examples to-day are complex structures on products of spheres of the form $S^1\times S^{2n-1}$ (and their complex submanifolds; see e.g. \cite{OV1}, \cite{OV2}) and some Oeljeklaus-Toma manifolds, which will be described  below. 

\subsection{Oeljeklaus-Toma manifolds}
We follow the original paper \cite{OeTo}. Fix a number field $K$ having $s>0$ real embeddings ang $2t>0$ complex embeddings. Let $\H$ be the complex upper half plane; then the ring of integers $\O(K)$ of $K$ acts on $\H^s\times \C^t$ by 
$$a\cdot (z_1,\dots,  z_{s+t})=\left(z_1+\sigma_1(a),\dots, z_{s+t}+\sigma_{s+t}(a)\right).$$

Next, the group of {\em totally positive units} $\O_K^{*, +}$ (i.e. units $u\in \O_K^*$ with positive value in all real embedings of $K$) also acts on $\H^s \times \C^t$ in a similar way by

$$u\cdot \left(z_1,\dots, z_{s+t}\right)=\left(\sigma_1(u)z_1,\dots, \sigma_{s+t}(u)z_{s+t}\right).$$
 If a subgroup $U\subset \O_K^{*, +}$  with $rank(U)=s$ is  such that its projection onto its first $s$ factors of its logarithmic embedding is a full lattice in $\R^s$ (such subgroups are called {\em admissible subgroups}) then combining the above action of $\O_K$ with the action of units in $U$ gives a co-compact, properly discontinous action of 
$U\ltimes \O_K$ on $\H^s\times \C^t;$ the resulting quotient will be denoted $X(K, U)$ and called an {\em Oeljeklaus-Toma manifold}. 

We recollect some facts about the manifolds $X(K, U);$ once again, we refer the original paper \cite{OeTo} for details and proofs.

\medskip

\theorem
 a) For any choice of the number field $K$ and of the admissible subgroup $U$, the manifold $X(K, U)$ is non-K\"ahler;

b) for $t=1, s>0$ and any choice of admissible $U$, the manifold $X(K, U)$ has an LCK structure;

c) for $s=1, t>1$ and any choice of admissible $U$, the manifold $X(K, U)$ has no LCK structure.

\hfill
\section{The results}
A classical theorem due   to L. Kronecker (in 1857) asserts that if a unit of some number field $K$  has the  same absolute value  in {\em all} the embeddings of $K$ must be a root of unity. 

Since then, many variations of this theme ({\em algebraic integres with specified restrictions on the absolute values of it Galois conjugates}) appeared.   As we shall see below, the search for LCK metrics on Oeljeklaus-Toma manifolds leads naturally to a (yet another) problem in this theme, namely the search for {\em units  of number fields with the same absolute value in all complex embeddings.}

This is due to the following:

\medskip

\lemma If an Oeljeklaus-Toma manifold $X=X(K, U)$ has an LCK metric, then its automorphy factors $\chi(u)$ for   $u\in U$  are given by:
\begin{eqnarray}\label{eqmod}
\chi(u)=\vert \sigma_{s+1}(u)\vert ^2=\dots=\vert \sigma_{s+t}(u)\vert ^2
\end{eqnarray}
for all $u\in U.$ In particular, for all $u\in U$ one has
$$\vert \sigma_{s+1}(u)\vert=\dots=\vert \sigma_{s+t}(u)\vert.$$

\hfill

\noindent {\bf Proof.}
Assume that $\omega$ is K\"ahler metric on $\H^s\times \C^t$ upon which $U\ltimes \O_K$ acts by homotheties. Then $\omega$ can be written as
{\small
$$\omega =
\sum_{i,j=\overline{1,s+t}} h_{i\overline{j}}dz_i\wedge d\overline{z}_j.
 $$
}
By an average argument, as in \cite{OeTo}, we can assume that all the coefficients of $\omega$ depend only on $z_1,\dots, z_s.$
Next, we infer that all $h_{i\overline{i}}, i=s+1,\dots, s+t$ are constant. Indeed, if this would not be the case, then
$$\frac{\partial h_{i\overline{i}}}{\partial z_{k}}\not =0 $$
  for some $1\leq k\leq s.$
But from the K\"ahlerianity assumption we get that also
$$\frac{\partial h_{k\overline{i}}}{\partial z_{i}}\not =0,$$
a contradiction with our assumption on $h_{k\overline{i}}.$
Now the conclusion on the automorphy factors follows at once.

\medskip

\section{The number-theoretical issues}
In this section, we state and prove the main number-theoretical ingredients needed for the proof of the main result. Since this section may be of interest for number-theorists (who may wish to skip the other sections), we recall the setup.

Fix  $K$ a number field with $s$ real mebedings and $2t$ complex embeddings; we label $\sigma_i, i=1,\dots, s+2t$ its embeedings, with the convention that the first $s$ ones are real, and for any $i=1,\dots, t$ one has $\sigma_{s+t+i}=\overline{\sigma}_{s+i}.$

We will introduce the following {\em ad-hoc} terminology, inspired by the equalities (\ref{eqmod}), which we consider suggestive for the geometrical context we are working in.
 
\medskip

\noindent  {\bf Definition.} A unit $u\in \O_K^*$ will be called 
{\em homothetical} if 
$$\vert\sigma_{s+1}\vert=\dots=\vert\sigma_{s+t}(u)\vert\not =1$$ and respectively {\em isometrical} if $$\vert\sigma_{s+1}\vert=\dots=\vert\sigma_{s+t}(u)\vert=1.$$

\medskip

We recall also the notion of {\em admissible subgroup} of units, needed for the construction of Oeljeklaus-Toma manifolds; it is a subgroup $U\subset \O_K^*$ whose projection onto its first $s$ factors in in its logarithmic embedding form a full lattice in $\R^s.$ 

Recall that for the construction of Oeljeklaus-Toma manifolds, we need admissible subgroups $U$ with  $rank(U)=s.$
The main scope of this section is to prove that for $1\leq s<t$ such subgroups do not exists.

\theorem
Let $X= X(K, U)$ be an Oeljeklaus-Toma manifold associated to a number field  $K$ with $s$ real embeddings and $2t$ complex embeddings and to an admissible subgroup $U\subset \O_K^{*, +}.$ If $1<s<t$ then $X$ has no lck metric.
\hfill

\noindent {\bf Proof.} 
The plan of proof is follows. Assuming $X$ has an LCK metric, we  prove first that all homothetical units have degree $<[K:\Q].$ This will imply that there is some proper subfield $L\subset K$ such that $U\subset \O_L^*.$ But this will force some of the complex embeddings of these units to have different absolute value, contradiction.

To prove the first assertion, let 
 $u\in U$ be a homothetical unit of maximal degree $deg(u)=[K:\Q]$,
 and let us label by $r_1,\dots, r_s$ its images under the real embeddings and respectively $z_1, \dots, z_{2t}$ its images under the complex embeddings (with the convention that $z_{t+i}=\overline{z}_{i}$ for all $i=1,\dots t$). Let us also  denote by $R$ the common value of $\vert z_k\vert, k=1,\dots, t.$

\noindent Recall that we have: 
\begin{eqnarray}\label{eqval}
z_1\overline{z}_1=\dots =z_t\overline{z}_t.
\end{eqnarray}
Let $K^{ncl}$ be the normal closure of $K;$ then for any $i=1,\dots, s$ there exists some $\sigma\in Gal(K^{ncl}\vert \Q)$ such that $\sigma(z_1)=r_i.$ Applying $\sigma$ to (\ref{eqval}) we get

\begin{eqnarray}\label{eq2}
r_i\sigma(\overline{z}_1)=\dots =\sigma(z_t)\sigma(\overline{z}_t.)
\end{eqnarray}
Since $s<t$, in the above equations (\ref{eq2}) we must have at least one factor of the form $\sigma(z_k)\sigma(\overline{z}_k)$ equal to some $z_\alpha z_\beta$  for some $k=2,\dots t$ and some $\alpha, \beta \in 1,\dots, 2t.$
Hence, we have

\begin{eqnarray}\label{dec}
r_i\sigma(\overline{z}_1)=z_\alpha z_\beta.
\end{eqnarray}
We consider the occuring possibilities.

\noindent {\bf Case 1.} There is some $\gamma(i)\in 1,\dots, 2t$ such that $\sigma(\overline{z}_1)=z_{\gamma(i)}.$
Taking absolute values, we get 
$r_i=R.$
But then, since $u$ was assumed to be of degree $$deg(u)=[K:\Q]$$ we have that all the $z_1,\dots, z_{2t} $ are distinct, so  by \cite{Boyd}, the minimal polynomial$f\in \Q[X]$  of $u$ is of the form
$f(X^{2t+1}).$ We get $2t+1\vert s+2t$; but this is absurd, as $1<s<t.$

So we are left with 

\noindent {\bf Case 2.} For all $i=1,\dots, s$, there exists $\f(i)\in 1,\dots, s$  and some $\alpha(i), \beta(i)\in 1,\dots 2t$ such that $r_ir_{\f(i)}=z_{\alpha(i)} z_{\beta(i)}$. Again, taking absolute values we get
$$r_ir_{\f(i)}=R^2.$$
Noticing that $\f$ is a bijection, (since by assumption $u$ was of maximal degree, so all the $r_i's$ are distinct), we this this implies 
$$\prod_{i=1}^s r_i=R^s.$$ But as $u$ is a totally positive unit, we have
$$\left(\prod_{i=1}^sr_i\right)R^{2t}=1,$$ hence we get $R=1$, again a contradiction, as $u$ was assumed to be a homothetical unit. We conclude that every homothetical unit has degree $<[K:\Q].$

\medskip

Next, let $L_1, \dots, L_M$ be the set of proper subfields of $K$ generated by the homothetical units (i.e. for each $i$ there is some homothetical unit $u_i$ such that $L_i=\Q(u_i)$) and for each $i=1,\dots, M$ let 
$$C_i=\{u\in U\vert u \text{=homothetical unit}, u\in L_i\}.$$
Let us also $Isom(U)$ for the subset of $U$ formed by the isometrical units; it is a proper subgroup of $U.$
As
$$\bigcup_{i=1}^M C_i=U\setminus Isom(U)$$
we see $U=<C_{i_0}>$ for some $i_0$ (where $<C_{i_0}>$ is the subgroup generated by $C_{i_0}$).
Hence, 
$$U\subset \O_{L_{i_0}}^*.$$
But then, at least two complex embeddings $\sigma_k, \sigma_l$ of $K$ lie over different real embeddings of $L_{i_0}.$ To see this, let $s'$ (respt $2t'$) be the number of real (respectively complex) embeddings of $L_{i_0}$ and let $l=[K:L_{i_0}].$
As $U$ is admissible, we must have 
$s'=s$ 
(cf \cite{OeTo}, Lemma 1.6) and the restriction of any two different real embeddings of $K$ to $L_{i_0}$ cannot coincide.   If at most one real embedding of $L$ extends to a complex embedding of $K$, then $s\geq l(s'-1).$  
We further get $$s\geq l(s-1)$$
If $s\geq 3$ we get $l=1$, a contradiction. For $s=2$, we get $l=2=s'$. 
But as $l=2$ and as any real emebdding of $L_{i_0}$ extends to at least one real embedding of $K,$ we get $s\geq 4,$ a contradiction again.

Hence  $\vert \sigma_k(u_{i_0})\vert \not = \vert \sigma_l(u_{i_0})\vert$, contradiction with the assumption on $u_{i_0}.$ Q.E.D.


{\small

\noindent {\sc Victor Vuletescu\\ University of Bucharest, Faculty of
Mathematics,
\\14
Academiei str., 70109 Bucharest, Romania.}\\
\tt vuli@fmi.unibuc.ro

}


\begin{thebibliography}{100}


\bibitem[Boyd]{Boyd} D. Boyd, {\em Irreducible polynomials with many roots of maximal modulus. }
Acta Arith. 68, No.1, 85--88 (1994).

\bibitem[Bel]{Bel}
F.A. Belgun, {\em On the metric structure of
non-K{\"a}hler complex surfaces}, Math. Ann. {\bf 317} (2000),
1--40.

\bibitem[Bru1]{Bru1}
M. Brunella, {\em Locally conformally Kähler metrics on certain non-Kählerian surfaces. }
Math. Ann. 346, No. 3, 629--639 (2010).

\bibitem[Bru2]{Bru2}
M. Brunella, {\em Locally conformally Kähler metrics on Kato surfaces.}
Nagoya Math. J. 202, 77--81 (2011).


\bibitem[GaOr]{GaOr} P.Gauduchon, L. Ornea:
{\em Locally conformal Kaehler metrics on Hopf surfaces}, Annales de l'Institut Fourier, 48 (1998), 1107--1127.

\bibitem[OeTo]{OeTo} K. Oeljeklaus, M. Toma:
{\em Non-K\"ahler compact complex manifolds associated to number fields. }
Ann. Inst. Fourier 55, No. 1, 161--171 (2005).

\bibitem[OV1]{OV1} L. Ornea, M. Verbitsky, {\em Structure theorem for
compact Vaisman manifolds}, Math. Res. Lett., {\bf 10} (2003), 799--805.

\bibitem[OV3]{OrVeFG} L. Ornea, M. Verbistky,
{\em Topology of locally conformal Kaehler manifolds with potential} 
International Mathematics Research Notices, 4 (2010), 117--126. 

\bibitem[OV2]{OV2} L. Ornea, M. Verbitsky, {\em A report on locally
conformally K\"ahler manifolds}, Contemporary Mathematics {\bf 542},
135-150, 2011.

\bibitem[Tri]{Tri} 
F. Tricerri, {\it Some examples of locally
conformal K{\"a}hler manifolds}, Rend.  Sem.  Mat.  Univ.
Politec.  Torino {\bf 40} (1982), 81--92.

\bibitem[Va]{vaisman}
I. Vaisman, {\em On locally and globally conformal
K\"ahler
manifolds}, Trans. Amer. Math. Soc. {\bf 262} (1980), 533--542.

\end{thebibliography}
\end{document}